# Generalized Reflection Coefficients and Inverse Factorization of Hermetian Block Toeplitz Matrix.

August 13, 2004

**Rami Kanhouche**[1]

### Abstract

A factorization of the inverse of a Hermetian positive definite matrix based on a diagonal by diagonal recurrence formulae permits the inversion of Block Toeplitz matrices, using only matrix-vector products, and with a complexity of $O\left(n_1^3 n_2^2\right)$, where $n_1$, is the block size, and $n_2$ is the block matrix block-size.

**Introduction.**

The techniques developed here are based on a generalization of the reflection coefficients or partial autocorrelation presented in [1], they appear when the well know Levinson recursion for inverting matrix Toeplitz [4][5] is extended to two families of bi orthogonal polynomials.

These generalized coefficients result to be two triangular arrays of numbers, which completely describe the structure of the related matrix. In the Toeplitz case, the coefficients along each diagonal are identical, so they are determined by two sequences of pairs of numbers. In the Toeplitz Hermetian case, the coefficients collapse to the classical reflection coefficients.

In this paper we consider Hermetian positive definite matrices. In this case, the reflection coefficients result to be a triangular array of pairs of numbers whose product is a positive number less than one. The families of related orthogonal polynomials determine two Cholesky factorizations of the inverse of the matrix.

These generated coefficients are computed recursively diagonal by diagonal. Setting some initial conditions on the principal diagonal, at each step of the recursion a whole diagonal is computed from the previous results.

In the Block Toeplitz case, the related two triangular arrays of the generalized coefficients keep the block structure. The levinson algorithm has been generalized to this case and the reflection coefficients results to be matrices which definition involves a square root of a matrix. Several methods for choosing a convenient square root have been developed, see Dégerine[2] and Delsartre [3].

Applying the algorithm obtained for p.d. matrices to the Block Toeplitz case, at each step of the recursion, instead of computing a whole diagonal, only $n_1$ coefficients are computed from the previous results. We find a factorization of the inverse with a complexity $O\left(n_1^3 n_2^2\right)$, where $n_2$ is the number of blocks in the matrix.

---

[1] PhD Student at Lab. CMLA, Ecole Normale Supérieure de Cachan, 61, avenue du Président Wilson, 94235 CACHAN Cedex, France. phone: +33-1-40112688, mobile: +33-6-62298219, fax: +33-1-47405901, e-mail: rami.kanhouche@cmla.ens-cachan.fr, kanram@wanadoo.fr.

**Notations.**

We consider the Hermetian positive definite (covariance) matrix,

$$R := \{r_{i,j}\}_{i,j:=0,1,\ldots n-1}, \qquad (1)$$

$r_{i,j} \in C,$

$n \in N.$

For each matrix $R$ we consider the sub matrices family:

$$R^{k,l} := \{b_{i,j}\}_{i,j:=0,\ldots l-k}, \quad b_{i,j} := r_{i+k,j+k} \quad k < l, \; k,l := 0\ldots n-1,$$

for each $R^{k,l}$ matrix we define the following block diagonal matrix:

$${}^0 R^{k,l} := \begin{bmatrix} 0_{k \times k} & & \\ & R^{k,l} & \\ & & 0_{(n-l) \times (n-l)} \end{bmatrix},$$

which contains $R^{k,l}$ diagonally positioned, with zero filled elsewhere.
the transport square matrix $J$ defined as:

$$J := \begin{bmatrix} 0 & \cdots & 0 & 1 \\ \vdots & \iddots & 1 & 0 \\ 0 & \iddots & \iddots & \vdots \\ 1 & 0 & \cdots & 0 \end{bmatrix}.$$

## 1 Reflection Coefficients and Choleskey Factorization of inverse for Hermetian Positive definite matrices.

For any Hermetian p.d matrix the inverse can be factorized into the form [1]

$$R^{-1} = \left(R^P\right)^* \left(D^P\right)^{-1} \left(R^P\right)^T = \left(R^Q\right)^* \left(D^Q\right)^{-1} \left(R^Q\right)^T, \qquad (2)$$

where both of $D^P = \left(R^P\right)^T (R)\left(R^P\right)^*$, and $D^Q = \left(R^Q\right)^T (R)\left(R^Q\right)^*,$ (3)

are diagonal, and $R^P$, $R^Q$ are lower triangular and upper triangular respectively, both formed from Vectors $p_{k,l}, q_{k,l}$, $k,l = 0,1,\ldots n-1$, $k \leq l$ according to

$$R^P = \left[p_{0,n-1}, p_{1,n-1}, \ldots, p_{n-1,n-1}\right], \qquad (4)$$

$$R^Q = \left[q_{0,1}, q_{0,2}, \ldots, q_{0,n-1}\right].$$

These vectors (or polynomials) are mutually orthogonal relatively to the matrix defined cross product:

$$\langle v_1, v_2 \rangle := v_1^T R \left(v_2\right)^*,$$

this can be expressed in the conditions:

$$\forall \, 0 \leq k, 0 \leq k', k < l, k' < l, l \leq n-1,$$

$$\langle p_{k,l}, p_{k',l} \rangle = 0 \; if \; k \neq k',$$

$$\langle p_{k,l}, p_{k',l} \rangle > 0 \; if \; k = k'.$$

and

$$\forall \, 0 \leq k, k < l, k < l', l' \leq n-1, l \leq n-1,$$
$$\langle q_{k,l}, q_{k,l'} \rangle = 0 \text{ if } l \neq l',$$
$$\langle q_{k,l}, q_{k,l'} \rangle > 0 \text{ if } l = l'.$$

The orthogonal polynomials $p_{k,l}, q_{k,l}$ are obtained for $R$ in a recursive manner. The positive definiteness of the matrix R permits us to establish a recurrence system in which the famous reflection coefficients are used. The Generalized reflection coefficients define the relation between a group of vectors $p_{k,k+d}, q_{k,k+d}$ and the next step group of vectors $p_{k,k+d+1}, q_{k,k+d+1}$, $d := 0,1\ldots n-2$, according to

$$p_{k,l} = p_{k,l-1} - a_{k,l} \, q_{k+1,l}, \tag{5}$$
$$q_{k,l} = q_{k+1,l} - a'_{k,l} \, p_{k,l-1}, \tag{6}$$

starting from the canonical basis vectors $p_{k,k} := e_k, q_{k,k} := e_k$, $k = 0,1,\ldots n-1$.

The Generalized reflection coefficients are complex numbers, obtained in each step according to:

$$a_{k,l} = \frac{p_{k,l-1}^T \, R \, e_l}{q_{k+1,l}^T \, R \, e_l}, \tag{7}$$

$$a'_{k,l} = \frac{q_{k+1,l}^T \, R \, e_k}{p_{k,l-1}^T \, R \, e_k}. \tag{8}$$

The definition of $v_{k,l} := q_{k,l}^T \, R \, e_l$, $v'_{k,l} := p_{k,l}^T \, R \, e_k$ installs the following recurrence system:

$$v_{k,l} = v_{k+1,l} (1 - a_{k,l} a'_{k,l}) \in \mathbf{R}^+, \tag{9}$$
$$v'_{k,l} = v'_{k,l-1} (1 - a_{k,l} a'_{k,l}) \in \mathbf{R}^+. \tag{10}$$

This permit us to avoid applying the product in the denominator in (7) and (8) at each step, while for numerators, the following hold:

$$p_{k,l-1}^T \, R \, e_l = \left( q_{k+1,l}^T \, R \, e_k \right)^*. \tag{11}$$

**2 Generalized Reflection Coefficients in Hermetian Block Toeplitz Case.**

In this section we will consider the more limiting case of Block Toeplitz matrix in which each block is not necessarily Toeplitz. Each matrix block is of size $n_1 \times n_1$, while the total block matrix is of size $n_1 n_2 \times n_1 n_2$. To help our approach, we notice that the covariance matrix elements $r_{i,j}$ in this context admit the following property.

*Property 1.*

$$r_{i,j} = r_{i \bmod n_1, \, j - i \sec n_1} \quad \text{if } i \leq j,$$
$$r_{i,j} = (r_{j,i})^* \quad \text{if } i > j.$$

*where $a \sec b := b \operatorname{int}(a/b)$, $\operatorname{int}(x)$ is the integral part of x, and $a \bmod b$ is equal to the remainder of division of a by b.*

*Proof.* From the T-B-T structure we got for $k := 0,1,\ldots n_1 - 1, l := 0,1\ldots n_1 n_2 - 1$:

Generalized Reflection Coefficients and Inverse Factorization of Hermetian Block Toeplitz Matrix

$$r_{k+tn_1, l+tn_1} = r_{k,l}, \quad t = 0,1,\cdots n_2 - 1 - \text{int}(l/n_1),$$

Setting $i := k + tn_1, j := l + tn_1$, we notice that $(i,j) \in [0, n_1 n_2 - 1] \times [0, n_1 n_2 - 1]$, and:

$$r_{i,j} = r_{i - tn_1, j - tn_1}.$$

By the definition of $i$, $i - tn_1 = i \bmod n_1$ and $tn_1 = i \sec n_1$.

The previous property is to simply, the direct mathematical representation, of the Hermetian Block Toeplitz Matrix case, the following properties; *Property 2* and *Property 3* will help us to consider the main recurrent system in the Block-Toeplitz case, segmented into recurrent sub matrices.

***Property 2.***
$$\left( p_{k,l}^{^0R^{sd}}, q_{k,l}^{^0R^{sd}}, a_{k,l}^{^0R^{sd}}, a'_{k,l}^{^0R^{sd}}, v_{k,l}^{^0R^{sd}}, v'_{k,l}^{^0R^{sd}} \right) = \left( p_{k,l}^{R}, q_{k,l}^{R}, a_{k,l}^{R}, a'_{k,l}^{R}, v_{k,l}^{R}, v'_{k,l}^{R} \right),$$

with condition that $s \leq k < l \leq d$.
The entities $p_{k,l}, q_{k,l}, a_{k,l}, a'_{k,l}, v_{k,l}, v'_{k,l}$ that corresponds to a given covariance matrix $M$, are noted respectively as $p_{k,l}^{M}, q_{k,l}^{M}, a_{k,l}^{M}, a'_{k,l}^{M}, v_{k,l}^{M}, v'_{k,l}^{M}$.

***Property 3.*** $R^{k,l} = R^{k \bmod n_1, l - k \sec n_1}$
*Proof.* for each of both matrices elements we apply *property 1*:
$$r_{i,j}^{R_{k,l}} = r_{k+i, k+j} = r_{(k+i) \bmod n_1, (k+j) - (k+i) \sec n_1}, i, j = 0,1\ldots, l-k,$$

$$r_{i,j}^{R^{k \bmod n_1, l-k \sec n_1}} = r_{k \bmod n_1 + i, k \bmod n_1 + j} = r_{(k \bmod n_1 + i) \bmod n_1, k \bmod n_1 + j - (k \bmod n_1 + i) \sec n_1},$$

from which our proof is completed.

Now we are ready to state the main result of this paper.

***Theorem 1.*** *In the case of Hermetian Block Toeplitz Matrix we got:*
$$a_{k,l} = a_{k \bmod n_1, l - k \sec n_1},$$
$$a'_{k,l} = a'_{k \bmod n_1, l - k \sec n_1},$$
$$p_{k,l} = U^{k \sec n_1} p_{k \bmod n_1, l - k \sec n_1},$$
$$q_{k,l} = U^{k \sec n_1} q_{k \bmod n_1, l - k \sec n_1},$$
$$v_{k,l} = v_{k \bmod n_1, l - k \sec n_1},$$
$$v'_{k,l} = v'_{k \bmod n_1, l - k \sec n_1},$$
*where $U$ is the $(n_1 n_2 \times n_1 n_2)$ shift matrix defined as,*

$$U = \begin{bmatrix} 0 & 0 & \cdots & 0 & 0 \\ 1 & 0 & \ddots & & 0 \\ 0 & \ddots & \ddots & 0 & \vdots \\ \vdots & \ddots & 1 & 0 & 0 \\ 0 & \cdots & 0 & 1 & 0 \end{bmatrix},$$

*with $U^0 = \mathbf{1}$.*

*Proof.* We start by noticing that
$\forall 0 \leq k < l \leq n_1 n_2 - 1 \, \exists s,d : s \leq k < l \leq d$, so that we can write

$$\left( p_{k,l}^{0_{R^{s\,d}}}, q_{k,l}^{0_{R^{s\,d}}}, a_{k,l}^{0_{R^{s\,d}}}, a'^{0_{R^{s\,d}}}_{k,l}, v_{k,l}^{0_{R^{s\,d}}}, v'^{0_{R^{s\,d}}}_{k,l} \right) =$$

$$\left( U_{d-s}^{s}\, p_{k-s,l-s}^{R^{s\,d}}, U_{d-s}^{s}\, q_{k-s,l-s}^{R^{s\,d}}, a_{k-s,l-s}^{R^{s\,d}}, a'^{R^{s\,d}}_{k-s,l-s}, v_{k-s,l-s}^{R^{s\,d}}, v'^{R^{s\,d}}_{k-s,l-s} \right),$$

where $U_x$ is a Matrix of size $(n_1 n_2) \times (x+1)$, defined as:

$$U_x := \begin{bmatrix} 0 & 0 & \cdots & 0 \\ 1 & 0 & \ddots & \vdots \\ 0 & 1 & \ddots & 0 \\ 0 & 0 & \ddots & 0 \\ \vdots & \vdots & \ddots & 1 \\ 0 & 0 & \cdots & 0 \end{bmatrix}, \quad U_x^0 := \begin{bmatrix} \mathbf{1}_{(x+1)\times(x+1)} \\ \mathbf{0} \end{bmatrix},$$

coupling this with *property 2*, we obtain:

$$\left( U_{d-s}^{s}\, p_{k-s,l-s}^{R^{s\,d}}, U_{d-s}^{s}\, q_{k-s,l-s}^{R^{s\,d}}, a_{k-s,l-s}^{R^{s\,d}}, a'^{R^{s\,d}}_{k-s,l-s}, v_{k-s,l-s}^{R^{s\,d}}, v'^{R^{s\,d}}_{k-s,l-s} \right) =$$

$$\left( p_{k,l}^{R}, q_{k,l}^{R}, a_{k,l}^{R}, a'^{R}_{k,l}, v_{k,l}^{R}, v'^{R}_{k,l} \right).$$

By setting $s = k \sec n_1$, admitting that $k - k \sec n_1 = k \bmod n_1$, and using property 3, we obtain:

$$\left( \begin{array}{l} U_{d-s}^{s}\, p_{k \bmod n_1, l - k \sec n_1}^{R^{0,d-s}}, U_{d-s}^{s}\, q_{k \bmod n_1, l - k \sec n_1}^{R^{0,d-s}}, \\ a_{k \bmod n_1, l - k \sec n_1}^{R^{0,d-s}}, a'^{R^{0,d-s}}_{k \bmod n_1, l - k \sec n_1}, v_{k \bmod n_1, l - k \sec n_1}^{R^{0,d-s}}, v'^{R^{0,d-s}}_{k \bmod n_1, l - k \sec n_1} \end{array} \right) =$$

$$\left( p_{k,l}^{R}, q_{k,l}^{R}, a_{k,l}^{R}, a'^{R}_{k,l}, v_{k,l}^{R}, v'^{R}_{k,l} \right).$$

Finally, our proof is completed by noticing that :

$$\left( \begin{array}{l} U_{d-s}^{s}\, p_{k \bmod n_1, l - k \sec n_1}^{R^{0,d-s}}, U_{d-s}^{s}\, q_{k \bmod n_1, l - k \sec n_1}^{R^{0,d-s}}, \\ a_{k \bmod n_1, l - k \sec n_1}^{R^{0,d-s}}, a'^{R^{0,d-s}}_{k \bmod n_1, l - k \sec n_1}, v_{k \bmod n_1, l - k \sec n_1}^{R^{0,d-s}}, v'^{R^{0,d-s}}_{k \bmod n_1, l - k \sec n_1} \end{array} \right) =$$

$$\left( \begin{array}{l} U^{s}\, p_{k \bmod n_1, l - k \sec n_1}^{0_{R^{0,d-s}}}, U^{s}\, q_{k \bmod n_1, l - k \sec n_1}^{0_{R^{0,d-s}}}, \\ a_{k \bmod n_1, l - k \sec n_1}^{0_{R^{0,d-s}}}, a'^{0_{R^{0,d-s}}}_{k \bmod n_1, l - k \sec n_1}, v_{k \bmod n_1, l - k \sec n_1}^{0_{R^{0,d-s}}}, v'^{0_{R^{0,d-s}}}_{k \bmod n_1, l - k \sec n_1} \end{array} \right) =$$

$$\left( \begin{array}{l} U^{s}\, p_{k \bmod n_1, l - k \sec n_1}^{R}, U^{s}\, q_{k \bmod n_1, l - k \sec n_1}^{R}, \\ a_{k \bmod n_1, l - k \sec n_1}^{R}, a'^{R}_{k \bmod n_1, l - k \sec n_1}, v_{k \bmod n_1, l - k \sec n_1}^{R}, v'^{R}_{k \bmod n_1, l - k \sec n_1} \end{array} \right).$$

We are also interested in the algorithmic aspect of theorem 1; more precisely the following theorem explains the resulting algorithmic details and calculus coast.

***Theorem 2.*** *In Hermetian Block Toeplitz case, the application of Recurrence relations (5)-(10) resume in the following Algorithm, with complexity $O\left(n_1^3 n_2^2\right)$:*
*First we rewrite Equations (5)-(10) as the subroutine,*

*Subroutine 1*

$$a_{k,l} = \frac{p_{k^0,l^-}^T R\, e_l}{v} \qquad a'_{k,l} = \frac{\hat{q}^T R\, e_k}{v'_{k^0,l^-}} \qquad (t2.1)$$

$$p_{k,l} = p_{k^0,l^-} - a_{k,l}\hat{q} \qquad q_{k,l} = \hat{q} - a'_{k,l}\, p_{k^0,l^-} \qquad (t2.2)$$

$$v_{k,l} = v\left(1 - a_{k,l} a'_{k,l}\right) \qquad v'_{k,l} = v'_{k^0,l^-}\left(1 - a_{k,l} a'_{k,l}\right) \qquad (t2.3)$$

and proceed as the following,

*Initialization:*
For $k = 0$ to $n_1 - 1$: {
$p_{k,k} = q_{k,k} = e_k$, $v_{k,k} = v'_{k,k} = r_{0,0}$
}

*Main routine:*
For $d_2 = 0$ to $n_2 - 1$ do : (loop 1)
{
   If $d_2$ not equal to zero:
   {
      For $d_1 = n_1 - 1$ to $0$ do : (Lower triangle loop)
      {
         For $u = 0$ to $n_1 - d_1 - 1$ do :
         {
         $k = u + d_1 \quad l = d_2 n_1 + u$
         $(k^0, l^-) = (k, l-1)$
         if $u$ equal to $n_1 - d_1 - 1$: {
         $(k^+, l^0) = ((k+1)\bmod n_1, l - (k+1)\sec n_1)$
         $\hat{q} = U^{n_1} q_{k^+, l^0}$, $v = v_{k^+, l^0}$        (t2.4)
         }
         else {
         $(k^+, l^0) = (k+1, l)$
         $\hat{q} = q_{k^+, l^0}$, $v = v_{k^+, l^0}$
         }
         *Apply Subroutine 1*
         }
      }
   }
   For $d_1 = 1$ to $n_1 - 1$ do : (Upper triangle loop)
   {
      For $u = 0$ to $n_1 - d_1 - 1$ do :
      {
         $k = u \quad l = d_2 n_1 + u + d_1$
         $(k^0, l^-) = (k, l-1)$

$$(k^+, l^0) = (k+1, l)$$
$$\hat{q} = q_{k^+,l^0}, v = v_{k^+,l^0}$$
*Apply Subroutine 1*
}
}
}

*Proof.* For any element $(i,j) \in T_d$,
$$T_d := \{(k,l) : k < l, (k,l) \in [0, n_1 d - 1] \times [0, n_1 d - 1]\}, d = 1, 2, \ldots n_2;$$
$(i \bmod d, j - i \sec d) \in P_d$,
$$P_d := \{(k,l) : k < l, (k,l) \in [0, n_1 - 1] \times [0, n_1 d - 1]\}$$
From Theorem 1, by obtaining $p_{k,l}, q_{k,l}, a_{k,l}, a'_{k,l}, v_{k,l}, v'_{k,l}$ for all subscripts $(k,l) \in P_d$ we obtain directly these values for $(k,l) \in T_d$.

The algorithm proceeds in calculating values for all subscripts contained in the following groups, with respective order:
$\Lambda_1, \Lambda_2, \ldots, \Lambda_{n_2}$, where $\Lambda_d := P_d \setminus P_{d-1}$, and $P_0 := \emptyset$.

In *Lower triangle loop* we obtain values for
$\forall (k,l) \in \Lambda_d^-$, where $\Lambda_d^- := \{(k,l) \in \Lambda_d : k \geq l \bmod n_1\}$.

To demonstrate this, first we notice that for any group defined as
$$\Lambda_d^{-,d_1} := \{(k,l) : (k,l) \in \Lambda_d^-, k - l \bmod n_1 = d_1\},$$
$$\Lambda_d^- = \bigcup_{d_1=0}^{n_1-1} \Lambda_d^{-,d_1},$$
we got, for $d_1 \in [0, n_1 - 1]$:
$$\Lambda_d^{-,d_1} = \{(k,l) : l \bmod n_1 \in [0, n_1 - 1 - d_1], k = l \bmod n_1 + d_1,\}.$$

By defining $u := l \bmod n_2$ we conclude that we did consider all $(k,l) \in \Lambda_d^-$.
Using the same logic we conclude that in *Upper triangle loop* we obtain values for
$\forall (k,l) \in \Lambda_d^+$, where $\Lambda_d^+ := \{(k,l) \in \Lambda_d : k < l \bmod n_1\}$. It is clear that $\Lambda_d = \Lambda_d^- \cup \Lambda_d^+$

In the *upper triangle loop* it is easy to verify that:
$\forall (k,l) \in \Lambda_d^+, \{(k, l-1), (k+1, l)\} \subset \Lambda_d$,

while in the *Lower triangle loop* this is not always true, more precisely for $k = n_1 - 1$, or/and $l \bmod n_1 = 0$. So we can write:
$\forall (k,l) \in \Lambda_d^-, k \neq n_1 - 1 \Rightarrow \{(k, l-1), (k+1, l)\} \subset P_d$.

Finally for the case of $k = n_1 - 1$, we apply theorem 1 to obtain the needed values, since: $\forall (k,l) \in \Lambda_d^-, \{(k, l-1), ((k+1) \bmod n_1, l \sec n_1)\} \subset P_d$.

For computing the complexity, we proceed as the following:
Each Entry in the routine will require operations of $O(l - k)$, by noting $c_1$ as the step constant, $opc$ as the total number of operation, the calculus cost take the form of:

Generalized Reflection Coefficients and Inverse Factorization of Hermetian Block Toeplitz Matrix

$$opc = \sum_{d_2=0}^{n_2-1} \sum_{d_1=1}^{n_1-1} \sum_{u=0}^{n_1-d_1-1} c_1(n_1 d_2 + d_1) + \sum_{d_2=1}^{n_2-1} \sum_{d_1=0}^{n_1-1} \sum_{u=0}^{n_1-d_1-1} c_1(n_1 d_2 - d_1),$$

$$= \sum_{d_1=1}^{n_1-1}(n_1-d_1-1)c_1(d_1) + \sum_{d_2=1}^{n_2-1}\sum_{d_1=1}^{n_1-1}(n_1-d_1-1)c_1(n_1 d_2 + d_1) + \sum_{d_2=1}^{n_2-1}\sum_{d_1=0}^{n_1-1}(n_1-d_1-1)c_1(n_1 d_2 - d_1),$$

$$= \sum_{d_1=1}^{n_1-1}(n_1-d_1-1)c_1(d_1) + \sum_{d_2=1}^{n_2-1}\sum_{d_1=1}^{n_1-1}\left[(n_1-d_1-1)c_1(n_1 d_2 + d_1) + (n_1-d_1-1)c_1(n_1 d_2 - d_1)\right] +$$

$$\sum_{d_2=1}^{n_2-1}(n_1-1)c_1(n_1 d_2)$$

$$= \sum_{d_1=1}^{n_1-1}(n_1-d_1-1)c_1(d_1) + \sum_{d_2=1}^{n_2-1} 2c_1(n_1 d_2) \sum_{d_1=1}^{n_1-1}(n_1-d_1-1) + \sum_{d_2=1}^{n_2-1}(n_1-1)c_1(n_1 d_2),$$

$$= (n_1-1)c_1 \sum_{d_1=1}^{n_1-1} d_1 - c_1 \sum_{d_1=1}^{n_1-1}(d_1)^2 + \sum_{d_2=1}^{n_2-1} 2c_1(n_1 d_2)(n_1-1)^2 - \sum_{d_2=1}^{n_2-1} 2c_1(n_1 d_2) \sum_{d_1=1}^{n_1-1} d_1 +$$

$$(n_1-1)c_1 n_1 \sum_{d_2=1}^{n_2-1} d_2,$$

$$= (n_1-1)c_1 \frac{1}{2}(n_1-1)n_1 - c_1 \frac{1}{6}(n_1-1)(n_1)(2n_1-1) +$$

$$2c_1 n_1 \frac{1}{2}(n_2-1)n_2 (n_1-1)^2 - 2c_1 n_1 \frac{1}{2}(n_2-1)n_2 \frac{1}{2}(n_1-1)n_1 + \quad (12)$$

$$(n_1-1)c_1 n_1 \frac{1}{2}(n_2-1)n_2.$$

From which we conclude the complexity order of $n_1^3 n_2^2$.

**Bezoutian Inversion Formulae.**

One can directly develop a bezoutian factorization of the inverse, depending on (2). Instead of that, and for the benefit of the reader, we will express our new formulae depending on the already existing work of [8].
In [9], the inverse of a Hermetian Block Toeplitz matrix can be expressed in the form of :

$$R^{-1} = (L_A)^H D(P_f^{-1}) L_A - (L_B)^H D(P_b^{-1}) L_B, \quad (13)$$

where both of $L_A$, and $L_B$ are block Toeplitz in the form :

$$L_A = \begin{bmatrix} I & A_1 & \cdots & A_{n_2-1} \\ 0 & I & \ddots & \vdots \\ \vdots & \ddots & \ddots & A_1 \\ 0 & \cdots & 0 & I \end{bmatrix}, L_B = \begin{bmatrix} 0 & B_{n_2-1} & \cdots & B_1 \\ 0 & 0 & \ddots & \vdots \\ \vdots & \ddots & \ddots & B_{n_2-1} \\ 0 & \cdots & 0 & 0 \end{bmatrix}. \quad (14)$$

While the entities $\{A_i\}, \{B_i\}, P_f, P_b$, are defined as the unique solution to the equations [6],[7]:

$$\begin{bmatrix} I & A_1 & \cdots & A_{n_2-1} \end{bmatrix} R = \begin{bmatrix} P_f & 0 & \cdots & 0 \end{bmatrix}, \quad (15)$$

$$\begin{bmatrix} B_{n_2-1} & \cdots & B_1 & I \end{bmatrix} R = \begin{bmatrix} 0 & \cdots & 0 & P_b \end{bmatrix}, \quad (16)$$
8

$D(M)$ is a block diagonal matrix containing M on its main block diagonal with zeros elsewhere.

The next theorem will help establish the correspondence between the WWR algorithm, and the Generalized Reflection Coefficients algorithm.

**Theorem.** *In the case of strongly regular Hermetian Block Toeplitz matrix, the inverse is equal to*

$$R^{-1} = (L_p)^H D(V'^{-1}) L_p - (L_q)^H D(V^{-1}) L_q, \qquad (13)$$

*where*

$$L_p = \begin{bmatrix} P_0^T & P_1^T & \cdots & P_{n_2-1}^T \\ 0 & P_0^T & \ddots & \vdots \\ \vdots & \ddots & \ddots & P_1^T \\ 0 & \cdots & 0 & P_0^T \end{bmatrix}, L_q = \begin{bmatrix} 0 & Q_{n_2-1}^T & \cdots & Q_1^T \\ 0 & 0 & \ddots & \vdots \\ \vdots & \ddots & \ddots & Q_{n_2-1}^T \\ 0 & \cdots & 0 & 0 \end{bmatrix}, \qquad (14)$$

*While the entities* $\{P_i\}, \{Q_i\}, V', V$ *are defined as:*

$$\begin{bmatrix} P_0 \\ P_1 \\ \vdots \\ P_{n_2-1} \end{bmatrix} := \begin{bmatrix} p_{0,n_1n_2-1} & p_{1,n_1n_2-1} & \cdots & p_{n_1-1,n_1n_2-1} \end{bmatrix}, \begin{bmatrix} Q_{n_2-1} \\ \vdots \\ Q_1 \\ Q_0 \end{bmatrix} := \begin{bmatrix} q_{0,n_1n_2-n_1} & \cdots & q_{0,n_1n_2-2} & q_{0,n_1n_2-1} \end{bmatrix},$$

*and* $V' := \begin{bmatrix} v'_{0,n_1n_2-1} & 0 & \cdots & 0 \\ 0 & v'_{1,n_1n_2-1} & \ddots & \vdots \\ \vdots & \ddots & \ddots & 0 \\ 0 & \cdots & 0 & v'_{n_1-1,n_1n_2-1} \end{bmatrix}, V := \begin{bmatrix} v_{0,n_1n_2-n_1} & 0 & \cdots & 0 \\ 0 & \ddots & \ddots & \vdots \\ \vdots & \ddots & v_{0,n_1n_2-2} & 0 \\ 0 & \cdots & 0 & v_{0,n_1n_2-1} \end{bmatrix}.$

*Proof.* It is clear that both of $P_0$, and $Q_0$ are lower, and upper triangular respectively. From (2), and (3) we can write that:

$$R \cdot \begin{bmatrix} P_0 \\ P_1 \\ \vdots \\ P_{n_2-1} \end{bmatrix}^* V'^{-1}(P_0)^T = \begin{bmatrix} I \\ 0 \\ \vdots \\ 0 \end{bmatrix}, R \cdot \begin{bmatrix} Q_{n_2-1} \\ \vdots \\ Q_1 \\ Q_0 \end{bmatrix}^* V^{-1}(Q_0)^T = \begin{bmatrix} 0 \\ \vdots \\ 0 \\ I \end{bmatrix}.$$

From which, since R is Hermetian, we can write that:

$$(P_0)^* V'^{-1} \begin{bmatrix} P_0^T & P_1^T & \cdots & P_{n_2-1}^T \end{bmatrix} R = \begin{bmatrix} I & 0 & \cdots & 0 \end{bmatrix},$$

$$(Q_0)^* V^{-1} \begin{bmatrix} Q_{n_2-1}^T & \cdots & Q_1^T & Q_0^T \end{bmatrix} R = \begin{bmatrix} 0 & \cdots & 0 & I \end{bmatrix}.$$

This leads us to the equalities of

$$(P_0)^* V'^{-1} P_0^T = (P_f)^{-1}, \qquad (14) \qquad\qquad (Q_0)^* V^{-1} Q_0^T = (P_b)^{-1}, \qquad (15)$$

and eventually to

$$A_i = (P_0^T)^{-1} P_i^T, \qquad (16) \qquad\qquad B_i = (Q_0^T)^{-1} Q_i^T, \qquad (17)$$

since $P_f = \left(P_0^T\right)^{-1} V' \left(P_0^*\right)^{-1}$, and $P_b = \left(Q_0^T\right)^{-1} V \left(Q_0^*\right)^{-1}$.

By replacing (14)-(17) into (13) , (14) we obtain the proof of (13).

**Conclusion.**

The previous presented results, explain clearly the structure of the Generalized Reflection coefficients, and their relevant orthogonal polynomials in the Hermetian Block-Toeplitz Case, while the reflection coefficients admit the same Block Toeplitz status, the polynomials $p_{k,l}$ $q_{k,l}$, show a block Toeplitz recurrence with an added shift between Blocks polynomials counterparts.

**Acknowledgment.**

I would like to thank G. Castro[2], for providing, so gracefully, advice, and valuable information, to this article.

**References.**

[1] G. Castro, "Coefficients de Réflexion Généralisée. Extension de Covariance Multidimensionnelles et Autres Applications," Ph.D. Thesis, Université de Paris-Sud Centre.

[2] S. Dégerine, "Sample partial autocorrelation function of a multivariate time series", J. Multivariate Anal. 50 (1994), pp. 294-313.

[3] PH. Delsarte, Y. Genin and Y. Kamp, "Schur parametrization of positive definite block-Toeplitz systems", SIAM J. Appl. Math., 36, 1 (1979), pp 34-45.

[4] Durbin J., (1960) The fitting of times series models. Rev. Inst. Statist., Vol.28 pp 233-244.

[5] Levinson N., (1947) The Wiener RMS (root-mean-square) error Criterion in filter design and prediction J. Math. Anal. Phys, Vol.25, pp 261-278.

[6] P. Wittle, "On the Fitting of Multivariate Autoregressions, and the Approximative Canonical Factorization of a Spectral Density Matrix", *Biometrica*, 50:129-134, 1963.

[7] R. A. Wiggins and E. A. Robinson, "Recursive Solution of the Multichannel Filtering Problem", *J. Geophys. Res.*, vol. 70, pp. 1885-1891, April 1965.

[8] I. C. Gohberg and G. Heining, "Inversion of finite Toeplitz matrices made of elements of a non-commutative algebra", *Rev. Roumaine Math. Pures Appl.*, XIX(5):623-663, 1974, (In Russian).

[2] Universidad Central de Venezuela, Laboratorio de Probabilidades, CEPE, Escuela de Matématicas UCV, Caracas Venezuela. tel: 58-2126051512, fax: 52-2126051190